\documentclass[]{interact}

\usepackage{epstopdf}
\usepackage[caption=false]{subfig}

\usepackage[all, pdf]{xy}

\usepackage[numbers,sort&compress]{natbib}
\bibpunct[, ]{[}{]}{,}{n}{,}{,}
\makeatletter
\def\NAT@def@citea{\def\@citea{\NAT@separator}}
\makeatother

\theoremstyle{plain}
\newtheorem{theorem}{Theorem}[section]
\newtheorem{lemma}[theorem]{Lemma}
\newtheorem{corollary}[theorem]{Corollary}
\newtheorem{proposition}[theorem]{Proposition}

\newtheorem{definition}[theorem]{Definition}
\newtheorem{example}[theorem]{Example}
\newtheorem{remark}[theorem]{Remark} 

\begin{document}

\title{Property (R) with its stability for analytic functional system and under compact perturbations}
\author{
\name{Sinan Qiu\textsuperscript{1}, Lining Jiang\textsuperscript{1,*}, Xiaomin Wei\textsuperscript{2}}\thanks{\textsuperscript{*}Corresponding author. Email: jianglining@bit.edu.cn}
\affil{1.School of Mathematics and Statistics,
  Beijing Institute of Technology,
  Beijing, People's Republic of China}
  \affil{2.School of Mathematical and Physical Sciences, Nanjing Tech University, Nanjing, People's Republic of China}}

\maketitle

\begin{abstract}
 Let $\mathcal{F}_T$ be the analytic functional system generated by a linear bounded operator $T$ on an infinite dimensional complex Hilbert space.
 This paper characterizes the stability of property (R), a variant of Weyl's theorem, for $\mathcal{F}_T$ and under compact perturbations in terms of the  spectrum on consistency in Weyl. In particular,  the stability of property (R) under compact perturbations is shown to imply the (R)-stability of $\mathcal{F}_T$.
\end{abstract}


\section{Introduction}

A classical result of H. Weyl asserts that the essential spectrum of a self-adjoint operator on a Hilbert space is invariant under compact self-adjoint perturbations. Equivalently, the intersection of the spectra of all such perturbations coincides with the essential spectrum of the original operator, that is, its spectrum excluding the isolated eigenvalues of finite multiplicity \cite{ca}. This result is fundamental to the spectral analysis of self-adjoint operators and has important applications in mathematical quantum mechanics. In particular, its extension to relatively compact perturbations provides a useful tool for determining the essential spectra of Schr$\ddot{\rm o}$dinger operators \cite{11}.

Weyl's theorem and its variants are collectively termed Weyl type theorems.
Property (R) (see \cite{14}), as a Weyl type theorem, has garnered significant interest among scholars in recent years (see \cite{14, 15, 16, ba}). Let $H$ be an infinite dimensional complex Hilbert space, and denote by $B(H)$ the algebra of linear bounded operators on $H$. For $T\in B(H)$, let Hol($\sigma(T)$) denote the set of all functions that are analytic on an open neighborhood of $\sigma(T)$. 
Ren characterized property (R) from an operator $T$ to its functional calculus $f(T)$, and described stability of property (R) under commuting power finite rank perturbations with the help of CI spectrum (see \cite{15}). 
Jia revealed stability of property (R) under compact perturbations (see \cite{ba}). 
Yang in \cite[Theorem 2.3]{16} further showed that:
If an operator $T$ possesses the stability of property (R) under compact perturbations,
so does $f(T)$ for any $f\in\rm{Hol}(\sigma(T))$.
 However, no scholar has established the relationship between the stability of property (R) under compact perturbations and the property (R) for functional calculus of $T$.
 In this paper, we introduce the notion of an analytic functional system to characterize the stability of property (R) under compact perturbations in terms of the collective spectral properties of the functional calculus.
 Denote by $\mathcal{F}_T$ the set \[\{f(T):f\in\rm{Hol}(\sigma(T))\},\] 
 which we call an analytic functional system generated by $T$. If $f(T)$ possesses property (R) for any $f(T)\in\mathcal{F}_T$, we call the system $\mathcal{F}_T$ is (R)-stable. This paper aims to characterize the stability of property (R) under compact perturbations through the (R)-stability of the system $\mathcal{F}_T$.

In \cite{aa}, Cao introduced the notion of consistency in Weyl. Building upon this concept, we in this paper investigate its corresponding spectral properties and firstly define the resulting induced spectrum, referred to as the CW spectrum.
With the help of the CW spectrum, we characterize (R)-stability of the system $\mathcal{F}_T$.
Finally, 
the CW spectrum is used to show the relationship between the stability of property (R) under compact perturbations and (R)-stability of $\mathcal{F}_T$.

Throughout this paper, ${\mathbb C}$ denotes the set of complex numbers, and ${\Bbb N}$ denotes the set of nonnegative integers.
Let $T\in B(H)$. The nullity of $T$, denoted by $n(T)$, is the dimension of its kernel $N(T)$.
The defect of $T$, denoted by $d(T)$, is the codimension of its range $R(T)$.
$\sigma(T)$ and $\sigma_p(T)$ are the spectrum and the point spectrum of $T$ respectively.
Let $SF_+(H)=\{T\in B(H):$ $n(T)<\infty$ and $R(T)$ is closed$\}$, $SF_-(H)=\{T\in B(H):$ $d(T)<\infty$$\}$ and $SF(H)=SF_+(H) \cup SF_-(H)$. 
$T$ is said to be upper (resp. lower) semi-Fredholm if $T\in SF_+(H)$ (resp. $T\in SF_-(H)$).
 If $T\in SF(H)$, $T$ is said to be semi-Fredholm. 
 For brevity, for $\lambda\in\Bbb C$, we write \(T-\lambda\) for \(T-\lambda I\).
The upper semi-Fredholm spectrum, the lower semi-Fredholm spectrum, and the semi-Fredholm spectrum of $T$ are denoted, respectively, by:
$\sigma_{SF_+}(T)=\{\lambda\in{\Bbb C}: T-\lambda \notin SF_+(H)\}$, 
$\sigma_{SF_-}(T)=\{\lambda\in{\Bbb C}: T-\lambda \notin SF_-(H)\}$,
$\sigma_{SF}(T)=\sigma_{SF_+}(T)\cap\sigma_{SF_-}(T)$.
Furthermore, $T$ is Fredholm if $T\in SF_+(H) \cap SF_-(H)$; its essential spectrum is defined as $\sigma_{e}(T)=\sigma_{SF_+}(T)\cup\sigma_{SF_-}(T)$.
If $T\in SF(H)$, the index of $T$ is defined as $\hbox{ind}T=n(T)-d(T)$.
$T$ is said to be bounded below if it is injective and has closed range.
The approximate point spectrum of $T$ is defined as $$\sigma_{a}(T)=\{\lambda\in{\Bbb C}: T-\lambda \hbox{ is not bounded below}\}.$$
$T$ is said to be Weyl if it is Fredholm with index zero; its Weyl spectrum is denoted by $\sigma_w(T)=\{\lambda\in{\Bbb C}: T-\lambda $ is not Weyl$\}$.

The ascent (resp. descent) of $T$ is defined as $\hbox{asc}(T)=\inf\{n\in{\Bbb N}: N(T^{n})=N(T^{n+1})\}$ (resp. $\hbox{des}(T)=\inf\{n\in{\Bbb N}: R(T^{n})=R(T^{n+1})\}$).
We define $\hbox{asc}(T)=\infty$ (resp. $\hbox{des}(T)=\infty$) if the ascent (resp. descent) is not finite.
$T$ is said to be upper semi-Browder if $T\in SF_+(H)$ with $\hbox{asc}(T)<\infty$.
The upper semi-Browder spectrum of $T$ is defined as $$\sigma_{ab}(T)=\{\lambda\in{\Bbb C}: T-\lambda \hbox{ is not upper semi-Browder}\}.$$
$T$ is said to be Browder if it is Fredholm of finite ascent and descent.
The Browder spectrum of $T$ is defined as $\sigma_{b}(T)=\{\lambda\in{\Bbb C}: T-\lambda \hbox{ is not Browder}\}$.

Set $\rho(T)={\Bbb C}\backslash\sigma(T)$ and $\rho_{\tau}(T)={\Bbb C}\backslash\sigma_{\tau}(T)$,
where $\tau\in\{ SF_+,\  SF_-, \ SF, \ e,\ a,\ w, \ ab, \\\ b\}$.
Let $\rho_{SF}^-(T)$ denote the set of $\lambda\in\mathbb C$ such that $T-\lambda $ is semi-Fredholm and  $\hbox{ind}(T-\lambda )<0$. 
$\sigma_{0}(T)$ is defined as the set of $\lambda\in\mathbb C$ for which $T-\lambda $ is Browder but not invertible. 

For a subset $E \subseteq{\Bbb C}$, we denote by $\hbox{iso}E$, $\partial E$, $\hbox{acc} E$ and $\hbox{int} E$ the set of isolated points, boundary points, accumulative points and interior points of $E$ respectively. Write $\Bbb D=\{\lambda\in{\Bbb C}: \lvert \lambda \rvert\leq1\}$ and $\Bbb T=\{\lambda\in{\Bbb C}: \lvert \lambda \rvert=1\}$.

\section{Consistency in Weyl and property (R) for linear bounded operators}

	An operator $T\in B(H)$ is said to possess property (R) (see \cite[Definition 2.3]{14})(abbrev. $T\in (R)$), 
	if $$\sigma_{a}(T)\backslash\sigma_{ab}(T)=\pi_{00}(T),$$ where $\pi_{00}(T)=\{\lambda\in{\rm iso\sigma(T)}: 0<n(T-\lambda )<\infty\}$.
This section introduces the spectral properties of consistency in Weyl and characterizes property (R) for linear bounded operators in terms of the spectrum induced by consistency in Weyl.

\begin{definition}\cite{aa}
An operator $T\in B(H)$ is said to be consistent in Weyl (abbrev. $T\in CW(H)$), if for any $S\in B(H)$, one of the following holds:
\begin{enumerate} 
\item Both $TS$ and $ST$ are Weyl;
\item Neither $TS$ nor $ST$ is Weyl.
\end{enumerate}
\end{definition}

The CW resolvent set of an operator $T$ is defined as $\rho_{{CW}}(T)=\{\lambda\in\mathbb{C} : T-\lambda \in CW(H)\}$. Its complement, denoted $\sigma_{CW}(T)=\mathbb{C}\setminus\rho_{_{CW}}(T)$, is called the CW spectrum of $T$.

\begin{lemma}\cite[Theorem 2.1]{aa} 
Let $T\in B(H)$. $T\in CW(H)$ if and only if one of the following holds:
\begin{enumerate} 
\item $T\in \Phi(H)$;
\item $R(T)$ is not closed;
\item $R(T)$ is closed and $n(T)=d(T)=\infty$.
\end{enumerate}
\end{lemma}

It follows from Lemma 2.2 that $T\notin CW(H)$ if and only if $T\in SF_+(H)\setminus SF_-(H)$ or $T\in SF_-(H)\setminus SF_+(H)$. 
In other words, $\sigma_{CW}(T)=[\rho_{{SF}_{+}}(T)\cap\sigma_{{SF}_{-}}(T)]\cup[\rho_{{SF}_{-}}(T)\cap\sigma_{{SF}_{+}}(T)]$. In addition, considering the perturbation theory of semi-Fredholm operators (see\cite{aj}), $\sigma_{CW}(T)\subseteq \Bbb C$ is open thus $\rho_{_{CW}}(T)$ is closed, and furthermore, $\hbox{iso} \sigma(T)\subseteq\hbox{int} \rho_{_{CW}}(T)$ and $\rho_{SF}(T)\subseteq\hbox{int} \rho_{_{CW}}(T)\cup\sigma_{CW}(T)$.

Based on the spectral properties of CW, property (R) can be described as follows.
\begin{proposition} 
\label{prop2.3}
Let $T\in B(H)$. Then $T\in(R)$ if and only if
 $\sigma_{b}(T)=[\rho_{_{CW}}(T)\cap\pi_{\infty}(T)]\cup\pi_{0}(T)\cup\hbox{acc}\sigma_a(T)\cup[\hbox{acc}\sigma(T)\cap\sigma_d(T)]$, where $\pi_{\infty}(T)=\{\lambda\in\Bbb C: n(T-\lambda )=\infty\}$, $\pi_{0}(T)=\{\lambda\in\sigma(T): n(T-\lambda )=0\}$ and $\sigma_d(T)=\{\lambda\in\Bbb C: R(T-\lambda )$ is not closed$\}$.
\end{proposition}

\begin{proof}
For necessity, it suffices to show $\sigma_{b}(T)\subseteq[\rho_{_{CW}}(T)\cap\pi_{\infty}(T)]\cup\pi_{0}(T)\cup\hbox{acc}\sigma_a(T)\cup[\hbox{acc}\sigma(T)\cap\sigma_d(T)]$. 
Suppose $\lambda_{0}\notin[\rho_{_{CW}}(T)\cap\pi_{\infty}(T)]
\cup \pi_{0}(T)\cup \hbox{acc}\sigma_{a}(T)\cup[\sigma_d(T)\cap\hbox{acc}\sigma(T)].$ Without loss of generality, assume
$\lambda_{0}\in\sigma(T)$, then $n(T-\lambda_0)>0$, which means $\lambda_{0}\in\sigma_a(T)$.
Now, claim $\lambda_0\notin\pi_{\infty}(T)$.
Indeed, assume on the contrary that $\lambda_0\in\pi_{\infty}(T)$, then $\lambda_0\notin\rho_{_{CW}}(T)$. According to the definition of $\sigma_{CW}(T)$, $T-\lambda_0\in SF_+(H)\setminus SF_-(H)$ or $T-\lambda_0\in SF_-(H)\setminus SF_+(H)$. 
Furthermore, it infers from $\lambda_0\notin\hbox{acc}\sigma_a(T)$ that $\hbox{ind}(T-\lambda_0)\leq 0$. 
Thus $T-\lambda_0I\in SF_+(H)$ and $d(T-\lambda_0)=\infty$, which contradicts with ``$\lambda_0\in\pi_{\infty}(T)$". The assertion is proved.
Hence $0<n(T-\lambda_0)<\infty$.
Based on these, if $\lambda_{0}\notin\hbox{acc}\sigma(T)$, then $\lambda_0\in\pi_{00}(T)$. 
Given that $T\in(R)$, $\lambda_0\notin\sigma_b(T)$. 
In the following, assume $\lambda_0\notin\sigma_d(T)$, which shows $T-\lambda_0\in SF_+(H)$. 
Meanwhile, note that $\lambda_0\notin\hbox{acc}\sigma_{a}(T)$, thus $\lambda_0\in\sigma_{a}(T)\backslash\sigma_{ab}(T)$. It further follows from $T\in(R)$ that $\lambda_0\notin\sigma_b(T)$.

As for sufficiency,
one can see that $\{[\sigma_{a}(T)\backslash\sigma_{ab}(T)]\cup\pi_{00}(T)\}\cap\{[\rho_{_{CW}}(T)\cap\pi_{\infty}(T)]\cup\pi_{0}(T)\cup\hbox{acc}\sigma_a(T)\cup[\hbox{acc}\sigma(T)\cap\sigma_d(T)]\}=\emptyset$, 
hence $[\sigma_a(T)\backslash\sigma_{ab}(T)\cup\pi_{00}(T)]\cap\sigma_b(T)=\emptyset$. In other words, $[\sigma_a(T)\backslash\sigma_{ab}(T)\cup\pi_{00}(T)]\subseteq\sigma_0(T)$. Thus $\sigma_a(T)\backslash\sigma_{ab}(T)=\pi_{00}(T)$, i.e., $T\in(R)$. 
\end{proof}

\begin{remark}
\label{re2.4}
\begin{enumerate}
\item It infers from Proposition 2.3 that $T\in(R)$ if $\sigma_{b}(T)=\partial\rho_{_{CW}}(T)\cup\hbox{acc}\sigma_{a}(T)\cup[\sigma_d(T)\cap\hbox{acc}\sigma(T)]$. However, $\sigma_{b}(T)$ may be not equivalent to $\partial\rho_{_{CW}}(T)\cup\hbox{acc}\sigma_{a}(T)\cup[\sigma_d(T)
\cap\hbox{acc}\sigma(T)]$ even if $T\in(R)$.
For instance, let $A, B\in B(\ell^2)$ be defined by: $A(x_1, x_2, x_3, \cdots)=(0, x_1, x_2, x_3, \cdots)$, $B(x_1, x_2, x_3, \cdots)=(0, x_2, 0, x_4, \cdots).$ Moreover, let $T=A\oplus B$.
One can check that $\sigma_{a}(T)=\sigma_{ab}(T)={\Bbb T}\cup\{0\}$ and $\pi_{00}(T)=\emptyset$, which shows $T\in(R)$. However, note that $\sigma_b(T)={\Bbb D}$ while $\partial\rho_{_{CW}}(T)\cup\hbox{acc}\sigma_{a}(T)\cup[\sigma_d(T)
\cap\hbox{acc}\sigma(T)]=\Bbb T$. Thus, $\sigma_{b}(T)\neq\partial\rho_{_{CW}}(T)\cup\hbox{acc}\sigma_{a}(T)\cup[\sigma_d(T)
\cap\hbox{acc}\sigma(T)]$.

\item If $T\in (R)$, then $\sigma_{b}(T)\subseteq\rho_{_{CW}}(T)
\cup\pi_{0}(T)\cup \hbox{acc}\sigma_{a}(T)\cup[\sigma_d(T)\cap\hbox{acc}\sigma(T)]$, not vice versa.
Indeed, let $T\in B(\ell^2)$ be defined by: $T(x_1, x_2, x_3, \cdots)=(0, 0, \frac{x_2}{2}, \frac{x_3}{3}, \cdots).$
One can see $\sigma_b(T)=\{0\}\subseteq\rho_{_{CW}}(T)\cup\pi_{0}(T)\cup \hbox{acc}\sigma_{a}(T)\cup[\sigma_d(T)\cap\hbox{acc}\sigma(T)]=\{0\}.$  However, note that $\sigma_a(T)=\sigma_{ab}(T)=\{0\}$ and $\pi_{00}(T)=\{0\}$, thus $T\notin(R)$.

\item Due to the definition of $\rho_{_{CW}}(T)$, $\rho_{_{CW}}(T)\cap\pi_{\infty}(T)=\rho_{_{CW}}(T)\cap\pi_{\infty}(T)\cap\sigma_{SF_-}(T)$. Thus $T\in(R)$ if and only if $\sigma_{b}(T)=[\rho_{_{CW}}(T)\cap\pi_{\infty}(T)\cap\sigma_{SF_-}(T)]
\cup \pi_{0}(T)\cup \hbox{acc}\sigma_{a}(T)\cup[\sigma_d(T)\cap\hbox{acc}\sigma(T)]$.

\item Note that $\pi_{\infty}(T)\subseteq\sigma_{SF_+}(T)$. Thus, $\sigma_{b}(T)=[\rho_{_{CW}}(T)\cap\sigma_{SF_+}(T)]
\cup\pi_{0}(T)\cup \hbox{acc}\sigma_{a}(T)\cup[\sigma_d(T)\cap\hbox{acc}\sigma(T)]$ if $T\in(R)$. However, the reverse doesn't hold. For instance,
$T\in B(\ell^2)$ is defined by $T(x_1, x_2, x_3, \cdots)=(0, 0, \frac{x_2}{2}, \frac{x_3}{3}, \cdots)$. Then 
$\sigma_{b}(T)=\{0\}=[\rho_{_{CW}}(T)\cap\sigma_{SF_+}(T)]
\cup\pi_{0}(T)\cup \hbox{acc}\sigma_{a}(T)\cup[\sigma_d(T)\cap\hbox{acc}\sigma(T)]$. Whereas, $\sigma_a(T)\backslash\sigma_{ab}(T)=\emptyset$ is not equal to $\pi_{00}(T)=\{0\}$, which shows $T\notin(R)$.

\item $T$ is isoloid (i.e., $\sigma_p(T)\subseteq\hbox{iso}\sigma(T)$) and belongs to property (R) if and only if $\sigma_{b}(T)=[\rho_{_{CW}}(T)\cap\pi_{\infty}(T)]\cup\hbox{acc}\sigma_{a}(T)\cup[\sigma_d(T)\cap\hbox{acc}\sigma(T)].$
\end{enumerate}
\end{remark}

 Remark 2.4 establishes a connection between $\partial\rho_{_{CW}}(T)$ and property (R). The boundary $\partial\rho_{_{CW}}(T)$ is further employed to characterize property (R) as follows.

\begin{corollary} 
	Let $T\in B(H)$. Then the following statements are equivalent.
	\begin{enumerate} 
		\item  $T\in(R)$;
		
		\item  $\sigma_b(T)=\partial\rho_{_{CW}}(T)\cup\pi_{0}(T)\cup\hbox{acc}\sigma_a(T)\cup
		[\sigma_d(T)\cap\hbox{acc}\sigma(T)]\cup\{\hbox{acc}[\rho_e(T)\cap\rho_a(T)]\cap\pi_{\infty}(T)\}$;
		
		\item  $\sigma_b(T)=\partial\rho_{_{CW}}(T)\cup\pi_{0}(T)\cup\hbox{acc}\sigma_a(T)\cup[\hbox{acc}\rho_e^-(T)\cap\sigma_{SF_+}(T)]\cup
		\pi_{\infty}(T)$, where $\rho_e^-(T)=\{\lambda\in\rho_e(T): \hbox{ind}(T-\lambda)<0\}$.
	\end{enumerate}
\end{corollary}

\begin{proof} 
	$(1)\Rightarrow(2)$: It infers from Proposition 2.3 that 
	$\sigma_{b}(T)=[\rho_{_{CW}}(T)\cap\pi_{\infty}(T)]
	\cup\pi_{0}(T)\cup \hbox{acc}\sigma_{a}(T)\cup[\sigma_d(T)\cap\hbox{acc}\sigma(T)].$
	Given that $\rho_{CW}(T)\cap\pi_{\infty}(T)=
	[\partial\rho_{_{CW}}(T)\cap\pi_{\infty}(T)]\cup[\hbox{int}\rho_{_{CW}}(T)\cap\pi_{\infty}(T)]
	\subseteq\partial\rho_{_{CW}}(T)\cup[\hbox{int}\rho_{_{CW}}(T)\cap\pi_{\infty}(T)],$
	and
	$\hbox{int}\rho_{_{CW}}(T)\cap\pi_{\infty}(T)\subseteq
	\hbox{acc}\sigma_a(T)\cup[\hbox{int}\rho_{_{CW}}(T)\cap\hbox{iso}\sigma_a(T)\cap\pi_{\infty}(T)]\subseteq\hbox{acc}\sigma_a(T)\cup[\hbox{acc}\{\rho_e(T)\cap\rho_a(T)]\cap\pi_{\infty}(T)\}.$
	Hence, $\sigma_b(T)\subseteq\partial\rho_{_{CW}}(T)\cup\pi_{0}(T)\cup\hbox{acc}\sigma_a(T)\cup
	[\sigma_d(T)\cap\hbox{acc}\sigma(T)]\cup\{\hbox{acc}[\rho_e(T)\cap\rho_a(T)]\cap
	\pi_{\infty}(T)\}$. The reverse is clear.
	
	$(1)\Rightarrow(3)$: Given that $\{\hbox{acc}[\rho_e(T)\cap\rho_a(T)]\cap
	\pi_{\infty}(T)\}\subseteq\pi_{\infty}(T)$ and 
	$\hbox{acc}\sigma(T)\subseteq\hbox{acc}\sigma_a(T)\cup\partial\rho_{_{CW}}(T)\cup\pi_{0}(T)\cup\hbox{acc}\rho_e^-(T)$, 
	it deduces that $\hbox{acc}\sigma(T)\cap\sigma_d(T)\subseteq\hbox{acc}\sigma_a(T)\cup\partial\rho_{_{CW}}(T)\cup\pi_{0}(T)\cup[\hbox{acc}\rho_e^-(T)\cap\sigma_d(T)]\subseteq\hbox{acc}\sigma_a(T)\cup\partial\rho_{_{CW}}(T)\cup\pi_{0}(T)\cup[\hbox{acc}\rho_e^-(T)\cap\sigma_{SF_+}(T)].$
	Therefore, $\sigma_b(T)\subseteq\partial\rho_{_{CW}}(T)\cup\hbox{acc}\sigma_a(T)
	\cup\pi_{0}(T)\cup[\hbox{acc}\rho_e^-(T)\cap\sigma_{SF_+}(T)]\cup\pi_{\infty}(T)$. The reverse is obvious.

The proofs of $(2)\Rightarrow(1)$ and $(3)\Rightarrow(1)$ are analogous to the proof of the sufficiency part of Proposition 2.3.
\end{proof}

Indeed, property (R) can also be described by the Weyl spectrum.
Based on Proposition 2.3 and the fact that $\sigma_w(T)\subseteq\sigma_{b}(T)$, 
$T\in(R)$ if and only if $\sigma_w(T)\subseteq[\rho_{_{CW}}(T)\cap\pi_{\infty}(T)]
\cup \pi_{0}(T)\cup \hbox{acc}\sigma_{a}(T)\cup[\sigma_d(T)\cap\hbox{acc}\sigma(T)].$
The relationship between $\sigma_w(T)$ and $[\rho_{_{CW}}(T)\cap\pi_{\infty}(T)]
\cup \pi_{0}(T)\cup \hbox{acc}\sigma_{a}(T)\cup[\sigma_d(T)\cap\hbox{acc}\sigma(T)]$
may potentially be a proper inclusion.
For instance, $A, B\in B(\ell^2)$ are defined by: $A(x_1, x_2, x_3, \cdots)=(0, x_1, x_2, x_3, \cdots)$ and $B(x_1, x_2, x_3, \cdots)=(x_2, x_3, x_4, \cdots).$ Put $T=A\oplus B$. 
Then $\sigma_{a}(T)=\sigma_{ab}(T)={\Bbb D}$ and $\pi_{00}(T)=\emptyset$, which shows $T\in(R)$. Now, $\sigma_w(T)={\Bbb T}$, furthermore, note that $[\rho_{_{CW}}(T)\cap\pi_{\infty}(T)]
\cup \pi_{0}(T)\cup \hbox{acc}\sigma_{a}(T)\cup[\sigma_d(T)\cap\hbox{acc}\sigma(T)]={\Bbb D}$, thus  $\sigma_w(T)\subsetneqq[\rho_{_{CW}}(T)\cap\pi_{\infty}(T)]
\cup \pi_{0}(T)\cup \hbox{acc}\sigma_{a}(T)\cup[\sigma_d(T)\cap\hbox{acc}\sigma(T)].$
In what follows, we discuss when $\sigma_w(T)$ equals $[\rho_{_{CW}}(T)\cap\pi_{\infty}(T)]
\cup \pi_{0}(T)\cup \hbox{acc}\sigma_{a}(T)\cup[\sigma_d(T)\cap\hbox{acc}\sigma(T)]$.

\begin{proposition} 
\label{prop2.7}
Let $T\in B(H)$. Then $T\in(R)$ if and only if
$\sigma_w(T)=[\rho_{_{CW}}(T)\cap\pi_{\infty}(T)]
\cup\pi_{0}(T)\cup[\hbox{acc}\sigma_{a}(T)\cap\pi^{\pm}(T)]
\cup[\sigma_d(T)\cap\hbox{acc}\sigma(T)]$, where $\pi^{\pm}(T)=\{\lambda\in\Bbb C: n(T-\lambda )\neq d(T-\lambda )\}$.
\end{proposition}

\begin{proof}
Suppose $T\in(R)$, it follows from Proposition 2.3 that $\sigma_{b}(T)=[\rho_{_{CW}}(T)\cap\pi_{\infty}(T)]
\cup \pi_{0}(T)\cup \hbox{acc}\sigma_{a}(T)\cup[\sigma_d(T)\cap\hbox{acc}\sigma(T)].$
Besides that, $\hbox{acc}\sigma_{a}(T)=[\hbox{acc}\sigma_{a}(T)\cap\pi^{\pm}(T)]\cup[\hbox{acc}\sigma_{a}(T)\cap\pi_{\infty}(T)\cap\sigma_{SF_-}(T)]\cup[\hbox{acc}\sigma_{a}(T)\cap\rho_w(T)]$,
where $[\hbox{acc}\sigma_{a}(T)\cap\pi_{\infty}(T)\cap\sigma_{SF_-}(T)]
\subseteq[\rho_{_{CW}}(T)\cap\pi_{\infty}(T)]$. Hence
$\sigma_w(T)\subseteq[\rho_{_{CW}}(T)\cap\pi_{\infty}(T)]
\cup\pi_{0}(T)\cup[\hbox{acc}\sigma_{a}(T)\cap\pi^{\pm}(T)]
\cup[\sigma_d(T)\cap\hbox{acc}\sigma(T)].$ The reverse is clear.

For sufficiency, one can see $\sigma_b(T)=\sigma_w(T)\cup[\sigma_b(T)\cap\rho_w(T)]=
\sigma_w(T)\cup[\hbox{acc}\sigma_a(T)\cap\rho_w(T)]
\subseteq\sigma_w(T)\cup[\hbox{acc}\sigma_a(T)\cap\pi^{\pm}(T)]\subseteq[\rho_{_{CW}}(T)\cap\pi_{\infty}(T)]
\cup \pi_{0}(T)\cup \hbox{acc}\sigma_{a}(T)\cup[\sigma_d(T)\cap\hbox{acc}\sigma(T)]$, and the reverse is clear. 
Accordingly,
$\sigma_b(T)=[\rho_{_{CW}}(T)\cap\pi_{\infty}(T)]
\cup\pi_{0}(T)\cup \hbox{acc}\sigma_{a}(T)\cup[\sigma_d(T)\cap\hbox{acc}\sigma(T)]$, which means $T\in(R)$ on account of Proposition 2.3. 
\end{proof}

 \begin{remark}
 \label{re2.8}
\begin{enumerate}

\item $T\in(R)$ if and only if $\sigma_w(T)=[\rho_{_{CW}}(T)\cap\pi_{\infty}(T)\cap\sigma_{SF_-}(T)]\cup\pi_{0}(T)\cup[\hbox{acc}\sigma_{a}(T)\cap\pi^{\pm}(T)]\cup[\sigma_d(T)\cap\hbox{acc}\sigma(T)]$.

\item $T$ is isoloid (i.e., $\sigma_p(T)\subseteq\hbox{iso}\sigma(T)$) and belongs to property (R)  if and only if $\sigma_w(T)=[\rho_{_{CW}}(T)\cap\pi_{\infty}(T)]\cup[\hbox{acc}\sigma_{a}(T)\cap\pi^{\pm}(T)]
\cup[\sigma_d(T)\cap\hbox{acc}\sigma(T)]$.

\end{enumerate}
\end{remark}

\begin{corollary} 
\label{co2.9}
Let $T\in B(H)$. Then the following statements are equivalent.
\begin{enumerate} 
	\item  $T\in(R)$;
	
	\item  $\sigma_w(T)=\partial\rho_{_{CW}}(T)
	\cup\pi_{0}(T)\cup[\hbox{acc}\sigma_{a}(T)\cap\pi_{\infty}(T)\cap\sigma_{SF_-}(T)]\cup[\hbox{acc}\sigma_{a}(T)\cap\pi^{\pm}(T)]
	\cup[\sigma_d(T)\cap\hbox{acc}\sigma(T)]\cup\{\hbox{acc}[\rho_e(T)\cap\rho_{a}(T)]\cap\pi_{\infty}(T)\}$;
	
	\item  $\sigma_w(T)=\partial\rho_{_{CW}}(T)
	\cup\pi_{0}(T)\cup[\hbox{acc}\sigma_{a}(T)\cap\pi^{\pm}(T)]\cup\pi_{\infty}(T)\cup
	[\hbox{acc}\rho_e^-(T)\cap\sigma_{SF_+}(T)]$.
\end{enumerate}
\end{corollary}

\begin{proof}
It suffices to show the proofs of $(1)\Rightarrow(2)$ and $(1)\Rightarrow(3)$.

For $(1)\Rightarrow(2)$, due to Corollary 2.5, $\sigma_w(T)\subseteq\partial\rho_{_{CW}}(T)
\cup\pi_{0}(T)\cup\hbox{acc}\sigma_{a}(T)
\cup[\sigma_d(T)\cap\hbox{acc}\sigma(T)]\cup\{\hbox{acc}[\rho_e(T)\cap\rho_{a}(T)]\cap\pi_{\infty}(T)\}$. Furthermore,
$\hbox{acc}\sigma_{a}(T)\subseteq[\hbox{acc}\sigma_{a}(T)\cap\pi_{\infty}(T)\cap\sigma_{SF_-}(T)]\cup[\hbox{acc}\sigma_{a}(T)\cap\pi^{\pm}(T)]\cup\rho_w(T).$
Thus $\sigma_w(T)\subseteq\partial\rho_{_{CW}}(T)
\cup\pi_{0}(T)\cup[\hbox{acc}\sigma_{a}(T)\cap\pi_{\infty}(T)\cap\sigma_{SF_-}(T)]\cup[\hbox{acc}\sigma_{a}(T)\cap\pi^{\pm}(T)]
\cup[\sigma_d(T)\cap\hbox{acc}\sigma(T)]\cup\{\hbox{acc}[\rho_e(T)\cap\rho_{a}(T)]\cap\pi_{\infty}(T)\}$. The reverse is obvious.

For $(1)\Rightarrow(3)$, from Corollary 2.5, $\sigma_w(T)\subseteq\partial\rho_{_{CW}}(T)\cup\pi_{0}(T)\cup\hbox{acc}\sigma_a(T)\cup[\hbox{acc}\rho_e^-(T)\cap\sigma_{SF_+}(T)]\cup
\pi_{\infty}(T)$. In addition, $\hbox{acc}\sigma_{a}(T)\subseteq[\hbox{acc}\sigma_{a}(T)\cap\pi_{\infty}(T)]\cup\pi_{\infty}(T)\cup\rho_w(T).$ So, $\sigma_w(T)\subseteq\partial\rho_{_{CW}}(T)
\cup\pi_{0}(T)\cup[\hbox{acc}\sigma_{a}(T)\cap\pi^{\pm}(T)]\cup\pi_{\infty}(T)\cup
[\hbox{acc}\rho_e^-(T)\cap\sigma_{SF_+}(T)]$.  The reverse is obvious.
\end{proof}

\section{(R)-stability of analytic functional system $\mathcal{F}_T$}

 This section aims to characterize the (R)-stability of the analytic functional system $\mathcal{F}_T$ with the help of the spectral properties of CW. We first present some examples.

\begin{example}
\label{ex3.1}
Let $A, B, C\in B(\ell^2)$ be defined by $A(x_1, x_2, x_3, \cdots)=(0, x_1, x_2, \cdots)$, $B(x_1, x_2, x_3, \cdots)=(0, x_1, \frac{x_2}{2}, \frac{x_3}{3}, \cdots)$, $C(x_1, x_2, x_3, \cdots)=(0, x_2, x_3, \cdots)$.
Let $T=(A+3)\oplus(B-1)\oplus(C+1)\oplus(C-3)$.
Then $\sigma_a(T)=\{\lambda\in\Bbb C:|\lambda-3|=1\}\cup\{-3,-2,-1,1\}$, $\sigma_{ab}(T)=\{\lambda\in\Bbb C:|\lambda-3|=1\}\cup\{-2,-1\}$ and $\pi_{00}(T)=\{-3,1\}$.
In addition, 
let $f(T)=T^2$. 
It follows from $\{-1,1\}\subseteq\hbox{iso}\sigma(T)$ that $1\in\hbox{iso}\sigma(f(T))$. 
Meanwhile, note that $n(f(T)-1)=n(T-1)+n(T+1)=1$, thus $1\in\pi_{00}(f(T))$. 
However, $1\notin\sigma_a(f(T))\backslash\sigma_{ab}(f(T))$ since $T+1\notin SF_+(\ell^2)$.
In other words, $\pi_{00}(f(T))\nsubseteq\sigma_a(f(T))\backslash\sigma_{ab}(f(T))$.
Besides that, 
given that $T-3$ is bounded below but not invertible, and $T+3$ is Browder and non-invertible,
we obtain that $9\in\sigma_a(f(T))\backslash\sigma_{ab}(f(T))$. However, $9\notin\pi_{00}(f(T))$ since 
$f(T)-9$ is not Browder. 
Hence $T\in(R)$ but $f(T)\notin(R)$.
\end{example}

\begin{example}
\label{ex3.2}
Let $A, B, C\in B(\ell^2)$ be defined by:
$A(x_1,x_2,x_3,\cdots)=(0,x_1,x_2,\cdots)$,
 $B(x_1, x_2, x_3, \cdots)=(0,0,\frac{x_2}{2},\frac{x_3}{3},\cdots)$,
$C(x_1,x_2,x_3,\cdots)=(x_1,0,0,\cdots),$
Let $T=( A+1)\oplus (B-1)\oplus C$.
 One can check that $\sigma_a(T^2)=\{re^{i\theta}: r=2(1+\cos\theta)\}\cup\{1\}=\sigma_{ab}(T^2)$, $\hbox{iso}\sigma(T^2)=\emptyset$. Thus $T^2\in(R)$.
However, $\sigma_a(T)=\{\lambda\in\Bbb C: |\lambda-1|=1\}\cup\{-1,1\}$,
$\sigma_{ab}(T)=\{\lambda\in\Bbb C: |\lambda-1|=1\}\cup\{-1\}$, $\pi_{00}(T)=\{-1\}$, which infers $\sigma_{a}(T)\setminus\sigma_{ab}(T)\nsubseteq \pi_{00}(T)$ and $\pi_{00}(T)\nsubseteq\sigma_{a}(T)\setminus\sigma_{ab}(T)$. Hence $T\notin(R)$.
\end{example}
As shown in the examples above, 
\begin{enumerate} 
\item $T\in(R)\nRightarrow $ $\mathcal{F}_T$ is (R)-stable;
\item $f(T)\in(R)$ for some $f\in$Hol$(\sigma(T))$ $\nRightarrow T\in(R)$.
  \end{enumerate} 
This naturally raises the question: Under what conditions is $\mathcal{F}_T$ (R)-stable? We begin with some preliminary results.

\begin{proposition} 
\label{FP}
Let $T\in B(H)$. 
\begin{enumerate} 
	\item If $\sigma_0(T)=\emptyset$, then $\mathcal{F}_T$ is (R)-stable if and only if $T\in(R)$;
	\item If $\sigma_0(T)\neq\emptyset$, then $\mathcal{F}_T$ is (R)-stable if and only if $\sigma_{ab}(T)=\sigma_b(T)$ and $\hbox{iso}\sigma(T)\backslash \pi_{\infty}(T)=\sigma_0(T)$.
\end{enumerate}
 \end{proposition}

\begin{proof} 
Suppose that the system $\mathcal{F}_T$ is (R)-stable.
It suffices to check that $\sigma_{ab}(T)=\sigma_b(T)$ and $\hbox{iso}\sigma(T)\backslash \pi_{\infty}(T)=\sigma_0(T)$ when $\sigma_0(T)\neq\emptyset$.
Suppose $\mu_0\in\sigma_0(T)$. Take arbitrarily $\lambda_0\notin\sigma_{ab}(T)$, and let 
$f_1(T)=(T-\mu_0 )(T-\lambda_0 )$. Then $0\notin\sigma_{ab}(f_1(T))$ and $n(f_1(T))=n(T-\lambda_0 )+n(T-\mu_0 )>0$, in other words, $0\in\sigma_a(f_1(T))\backslash\sigma_{ab}(f_1(T))$. Since $\mathcal{F}_T$ is (R)-stable, $f_1(T)$ is Browder, which shows $T-\lambda_0 $ is Browder. 
Thus $\sigma_{ab}(T)=\sigma_b(T)$.
 Furthermore, take arbitrarily $\nu_0\in\hbox{iso}\sigma(T)\backslash \pi_{\infty}(T)$, and let 
$f_2(T)=(T-\mu_0 )(T-\nu_0 )$. Then $0\in\pi_{00}(f_2(T))$. Given that $\mathcal{F}_T$ is (R)-stable,
thus $f_2(T)$ is Browder, from which one can see that $\nu_0\notin\sigma_b(T)$. Hence $\hbox{iso}\sigma(T)\backslash \pi_{\infty}(T)=\sigma_0(T)$.

The sufficiency is discussed in the following.
Assume that $\sigma_0(T)=\emptyset$. Since $T\in(R)$, $\sigma_a(T)=\sigma_{ab}(T)$ and $\pi_{00}(T)=\emptyset$. Take arbitrarily $f\in\hbox{Hol}(\sigma(T))$. Considering the spectral mapping theorems of $\sigma_a(\cdot)$ and $\sigma_{ab}(\cdot)$,  $\sigma_a(f(T))=f(\sigma_a(T))=f(\sigma_{ab}(T))=\sigma_{ab}(f(T))$. 
Besides that,
note that $\pi_{00}(f(T))\subseteq f(\pi_{00}(T))=\emptyset$, hence $\pi_{00}(f(T))=\emptyset$, which means $f(T)\in(R)$. Due to the arbitrariness of $f$, $\mathcal{F}_T$ is (R)-stable.
In what follows, assume $\sigma_0(T)\neq\emptyset$, then $\sigma_{ab}(T)=\sigma_b(T)$ and $\hbox{iso}\sigma(T)\backslash \pi_{\infty}(T)=\sigma_0(T)$. For any $f\in\hbox{Hol}(\sigma(T))$, $\sigma_{ab}(f(T))=f(\sigma_{ab}(T))=f(\sigma_b(T))=\sigma_b(f(T))$, 
thus $\sigma_a(f(T))\backslash\sigma_{ab}(f(T))\subseteq\pi_{00}(f(T))$. 
Besides that, for any $\mu\in\pi_{00}(f(T))$, 
suppose $$f(T)-\mu =a(T-\lambda_1)^{n_1}(T-\lambda_2)^{n_2}\cdots(T-\lambda_t)^{n_t}r(T),$$ where  $r(T)$ is invertible, $a\neq0$ and $\lambda_i\neq\lambda_j$($i\neq j$).
 It deduces that $\lambda_i\in\hbox{iso}\sigma(T)\cup\rho(T)$$(i=1,2,\cdots,t)$ and $n(T-\lambda_i )<\infty$, from which we obtain that $\lambda_i\in\{\lambda\in\hbox{iso}\sigma(T):n(T-\lambda )<\infty\}\cup\rho(T)\subseteq\rho_b(T)$.
 Hence $f(T)-\mu $ is Browder, which means $\mu\in\sigma_a(f(T))\backslash\sigma_{ab}(f(T))$.
 So, $\pi_{00}(f(T))\subseteq\sigma_a(f(T))\backslash\sigma_{ab}(f(T))$. 
 In summary, $f(T)\in(R)$. By again exploiting the arbitrariness of $f$, $\mathcal{F}_T$ is (R)-stable.
\end{proof}

The preceding proposition characterizes the (R)-stability of the analytic functional system. 
To streamline Proposition 3.3, we invoke the spectral properties of CW to replace the two spectral equalities $\sigma_{ab}(T)=\sigma_b(T)$ and $\hbox{iso}\sigma(T)\backslash \pi_{\infty}(T)=\sigma_0(T)$ with a single spectral condition in the following.

\begin{lemma} 
\label{lem3.4}
 Let $T\in B(H)$. Then the following statements are equivalent.
 \begin{enumerate} 
\item  $\sigma_{ab}(T)=\sigma_b(T)$ and $\hbox{iso}\sigma(T)\backslash \pi_{\infty}(T)=\sigma_0(T)$;

\item  $\sigma_b(T)=\partial\rho_{_{CW}}(T)\cup\hbox{acc}\sigma_{SF_+}(T)\cup[\rho_{SF_+}(T)\cap\{\lambda\in\Bbb C: \hbox{asc}(T-\lambda )=\infty\}]\cup\{\hbox{acc}[\rho_e(T)\cap\sigma(T)]\cap\sigma_d(T)\}\cup\pi_{\infty}(T)$;

\item  $\sigma_b(T)=[\rho_{_{CW}}(T)\cap\pi_{\infty}(T)]\cup\hbox{acc}\sigma_a(T)$;

\item  $\sigma_b(T)=\partial\rho_{_{CW}}(T)\cup\hbox{acc}\sigma_a(T)\cup\{\hbox{acc}[\rho_e(T)\cap\rho_a(T)]\cap\pi_{\infty}(T)\}$;

\item  $\sigma_b(T)=\partial\rho_{_{CW}}(T)\cup\hbox{acc}\sigma_a(T)\cup\pi_{\infty}(T)$.
\end{enumerate}
 \end{lemma}

\begin{proof} 
 $(1)\Rightarrow(2)$: It suffices to show $\sigma_b(T)\subseteq\partial\rho_{_{CW}}(T)\cup\hbox{acc}\sigma_{SF_+}(T)\cup[\rho_{SF_+}(T)\cap\{\lambda\in\Bbb C: \hbox{asc}(T-\lambda )=\infty\}]\cup[\hbox{acc}\{\rho_e(T)\cap\sigma(T)\}\cap\sigma_d(T)]\cup\pi_{\infty}(T)$.
Take arbitrarily $\lambda_0\notin\partial\rho_{_{CW}}(T)\cup\hbox{acc}\sigma_{SF_+}(T)\cup[\rho_{SF_+}(T)\cap\{\lambda\in\Bbb C: \hbox{asc}(T-\lambda )=\infty\}]\cup\hbox{acc}[\{\rho_e(T)\cap\sigma(T)\}\cap\sigma_d(T)]\cup\pi_{\infty}(T)$, then $n(T-\lambda_0 )<\infty$. 
We may assume that $\lambda_0\in\sigma(T)$. Now claim that $\lambda_0\in\hbox{int}\rho_{_{CW}}(T)$.
Indeed, assume on the contrary that $\lambda_0\notin\hbox{int}\rho_{_{CW}}(T)$. Since $\lambda_0\notin\partial\rho_{_{CW}}(T)$,
 $\lambda_0\in\sigma_{_{CW}}(T)$, 
which means
$T-\lambda_0 \in SF_+(H)\setminus SF_-(H)$ or $T-\lambda_0 \in SF_-(H)\setminus SF_+(H)$. 
Note that $\lambda_0\notin\hbox{acc}\sigma_{SF_+}(T)$, thus $T-\lambda_0 \in SF_+(H)$ and $d(T-\lambda_0 )=\infty$.
In addition, 
considering that $\lambda_0\notin\rho_{SF_+}(T)\cap\{\lambda\in\Bbb C: \hbox{asc}(T-\lambda )=\infty\}$, 
one can see $\lambda_0\notin\sigma_{ab}(T)$. Combining with $\sigma_{ab}(T)=\sigma_b(T)$,
$T-\lambda_0 $ is Browder, which contradicts ``$d(T-\lambda_0 )=\infty$". The assertion is proved.
Note that $\lambda_0\in\hbox{int}\rho_{_{CW}}(T)$ and $\lambda_0\notin\hbox{acc}\sigma_{SF_+}(T)$. Thus, there exists $\varepsilon>0$ such that $T-\lambda $ is Fredholm whenever $0<\lvert\lambda-\lambda_0\rvert<\varepsilon$, which deduces $\lambda_0\in\hbox{acc}[\rho_e(T)\cap\sigma(T)]\cup\hbox{iso}\sigma(T)$. In the following, we claim $\lambda_0\in\hbox{iso}\sigma(T)$. 
 Indeed, assume that $\lambda_0\in\hbox{acc}[\rho_e(T)\cap\sigma(T)]$. Since $\lambda_0\notin\hbox{acc}[\rho_e(T)\cap\sigma(T)]\cap\sigma_d(T)$, 
it infers that $R(T-\lambda_0 )$ is closed, 
from which one can see that $T-\lambda_0\in SF_+(H)$. Furthermore, given that $\lambda_0\notin\rho_{SF_+}(T)\cap\{\lambda\in\Bbb C: \hbox{asc}(T-\lambda )=\infty\}$, then $\lambda_0\notin\sigma_{ab}(T)=\sigma_b(T)$, which contradicts ``$\lambda_0\in\hbox{acc}[\rho_e(T)\cap\sigma(T)]\subseteq\hbox{acc}\sigma(T)$". The claim is proved.
Therefore, $\lambda_0\in\hbox{iso}\sigma(T)\backslash \pi_{\infty}(T)=\sigma_0(T)$, which shows $\lambda_0\notin\sigma_b(T)$.

 $(2)\Rightarrow(1)$: 
Since $\{\rho_{ab}(T)\cup[\hbox{iso}\sigma(T)\backslash \pi_{\infty}(T)]\}\cap\big\{\partial\rho_{_{CW}}(T)\cup\hbox{acc}\sigma_{SF_+}(T)\cup[\rho_{SF_+}(T)\cap\{\lambda\in\Bbb C: \hbox{asc}(T-\lambda )=\infty\}]\cup\{\hbox{acc}[\rho_e(T)\cap\sigma(T)]\cap\sigma_d(T)\}\cup \pi_{\infty}(T)\big\}=\emptyset$, it deduces that $\{\rho_{ab}(T)\cup[\hbox{iso}\sigma(T)\backslash \pi_{\infty}(T)]\}\subseteq\rho_b(T)$. Moreover, given that $\sigma_{ab}(T)\subseteq\sigma_b(T)$ and $\sigma_0(T)\subseteq\hbox{iso}\sigma(T)\backslash \pi_{\infty}(T)$, thus $\sigma_{ab}(T)=\sigma_b(T)$ and $\hbox{iso}\sigma(T)\backslash \pi_{\infty}(T)=\sigma_0(T)$.

 $(1)\Rightarrow(3)$: It follows from $(1)$ that $T\in(R)$. Combining with Proposition 2.3, 
$\sigma_b(T)=[\rho_{_{CW}}(T)\cap\pi_{\infty}(T)]\cup\hbox{acc}\sigma_a(T)\cup[\hbox{acc}\sigma(T)\cap\sigma_d(T)]\cup\pi_{0}(T)$.
Meanwhile, note that $\sigma_{ab}(T)=\sigma_b(T)$, then $\sigma_a(T)=\sigma(T)$, from which one can see $\hbox{acc}\sigma(T)\cap\sigma_d(T)\subseteq\hbox{acc}\sigma_a(T)$.
Besides that, since $\hbox{iso}\sigma(T)\backslash \pi_{\infty}(T)=\sigma_0(T)$, $\pi_{0}(T)\subseteq\hbox{acc}\sigma(T)$. Furthermore, due to $\sigma_a(T)=\sigma(T)$, $\pi_{0}(T)\subseteq\hbox{acc}\sigma(T)\cap\sigma_d(T)$. Thus $\sigma_b(T)\subseteq[\rho_{_{CW}}(T)\cap\pi_{\infty}(T)]\cup\hbox{acc}\sigma_a(T)$, and the reverse is clear. 

$(1)\Rightarrow(4)$: According to $T\in(R)$ and Corollary 2.5, $\sigma_b(T)=\partial\rho_{_{CW}}(T)\cup\pi_{0}(T)\cup\hbox{acc}\sigma_a(T)\cup
[\sigma_d(T)\cap\hbox{acc}\sigma(T)]\cup[\hbox{acc}(\rho_e(T)\cap\rho_a(T))\cap\pi_{\infty}(T)]$.
In addition, the process of $(1)\Rightarrow(3)$ shows $\pi_{0}(T)\subseteq\hbox{acc}\sigma(T)\cap\sigma_d(T)$ and $\hbox{acc}\sigma(T)\cap\sigma_d(T)\subseteq\hbox{acc}\sigma_a(T)$. Thus $(4)$ holds.

$(1)\Rightarrow(5)$: Analogous to $(1)\Rightarrow(3)$, it follows from Corollary 2.6 that $(5)$ holds.

The proofs of $(3)\Rightarrow(1)$, $(4)\Rightarrow(1)$ and $(5)\Rightarrow(1)$ are similar to that of $(2)\Rightarrow(1)$.
\end{proof}

Based on Proposition \ref{FP} and Lemma \ref{lem3.4}, the (R)-stability of $\mathcal{F}_T$ is characterized through the spectral properties of CW as follows.

\begin{theorem} 
\label{FTHM}
Let $T\in B(H)$. 
\begin{enumerate} 
	\item If $\sigma_0(T)=\emptyset$, then $\mathcal{F}_T$ is (R)-stable if and only if $T\in(R)$;
	\item If $\sigma_0(T)\neq\emptyset$, then $\mathcal{F}_T$ is (R)-stable if and only if $\sigma_b(T)=[\rho_{_{CW}}(T)\cap\pi_{\infty}(T)]\cup\hbox{acc}\sigma_a(T)$.
\end{enumerate}
 \end{theorem}

Indeed, Theorem \ref{FTHM} characterize the (R)-stability of $\mathcal{F}_T$ in terms of the relationship between $\rho_{_{CW}}(T)$ and $\sigma_b(T)$. 
By analogy with the characterization of property (R) for an operator $T$ in terms of $\rho_{_{CW}}(T)$ and $\sigma_w(T)$ in Section 2, we now proceed to investigate the (R)-stability of the system $\mathcal{F}_T$ using the relationship between $\rho_{_{CW}}(T)$ and $\sigma_w(T)$. To this end, the following preliminaries are required.

\begin{lemma} 
\label{lem3.11}
Let $T\in B(H)$. Then the following statements are equivalent.
\begin{enumerate} 
\item $\sigma_{ab}(T)=\sigma_b(T)$ and $\hbox{iso}\sigma(T)\backslash \pi_{\infty}(T)=\sigma_0(T)$;

\item $\sigma_w(T)=[\rho_{_{CW}}(T)\cap\pi_{\infty}(T)]\cup[\hbox{acc}\sigma_a(T)\cap\pi^{\pm}(T)]$;

\item $\sigma_w(T)=\partial\rho_{_{CW}}(T)\cup[\hbox{acc}\sigma_a(T)\cap\pi_{\infty}(T)\cap\sigma_{SF_-}(T)]\cup[\hbox{acc}\sigma_a(T)\cap\pi^{\pm}(T)]\cup\{\hbox{acc}[\rho_e(T)\cap\rho_a(T)]\cap\pi_{\infty}(T)\}$;

\item $\sigma_w(T)=\partial\rho_{_{CW}}(T)\cup[\hbox{acc}\sigma_a(T)\cap\pi^{\pm}(T)]\cup\pi_{\infty}(T)$.
\end{enumerate} 
\end{lemma} 

\begin{proof} 
Suppose $(1)$ holds. Considering Proposition 2.6 and Corollary 2.8, one can see
$\pi_{0}(T)\subseteq[\hbox{acc}\sigma(T)\cap\sigma_d(T)]$ and $\sigma_a(T)=\sigma(T)$, and $[\hbox{acc}\sigma(T)\cap\sigma_d(T)]=[\hbox{acc}\sigma_a(T)\cap\sigma_d(T)\cap\pi^{\pm}(T)]\cup[\hbox{acc}\sigma_a(T)\cap\sigma_d(T)\cap\pi_{\infty}(T)\cap\sigma_{SF_-}(T)]$. Analogous to Lemma \ref{lem3.4}, $(1)\Rightarrow(2)$ and $(1)\Rightarrow(3)$ can be proved. Besides that, $(2)\Rightarrow(1)$, $(3)\Rightarrow(1)$ and $(4)\Rightarrow(1)$ are also similar to Lemma \ref{lem3.4}. 

In what follows, it suffices to show $(1)\Rightarrow(4)$.
Given that $T\in(R)$, it deduces from Corollary 2.8 that $\sigma_w(T)=\partial\rho_{_{CW}}(T)\cup\pi_{0}(T)\cup[\hbox{acc}\sigma_a(T)\cap\pi^{\pm}(T)]\cup[\hbox{acc}\rho_e^-(T)\cap\sigma_{SF_+}(T)]\cup\pi_{\infty}(T)$. Since $\hbox{acc}[\rho_e^-(T)\cap\sigma_{SF_+}(T)]\subseteq[\hbox{acc}\rho_e^-(T)\cap\sigma_d(T)]\cup[\hbox{acc}\rho_e^-(T)\cap\pi_{\infty}(T)]\subseteq[\hbox{acc}\sigma(T)\cap\sigma_d(T)]\cup\pi_{\infty}(T)$, one can obtain that $[\hbox{acc}\rho_e^-(T)\cap\sigma_{SF_+}(T)]\cup\pi_{0}(T)\subseteq[\hbox{acc}\sigma_a(T)\cap\pi^{\pm}(T)]\cup\pi_{\infty}(T)$. Thus $(4)$ holds. 
\end{proof}

According to Lemma \ref{lem3.11}, the (R)-stability of $\mathcal{F}_T$ is characterized with the help of $\rho_{_{CW}}(T)$ and $\sigma_w(T)$.

\begin{theorem} 
\label{thm3.12}
Let $T\in B(H)$. 
\begin{enumerate} 
\item If $\sigma_0(T)=\emptyset$, then $\mathcal{F}_T$ is (R)-stable if and only if $T\in(R)$;
\item If $\sigma_0(T)\neq\emptyset$, then $\mathcal{F}_T$ is (R)-stable if and only if $\sigma_w(T)=[\rho_{_{CW}}(T)\cap\pi_{\infty}(T)]\cup[\hbox{acc}\sigma_a(T)\cap\pi^{\pm}(T)]$.
\end{enumerate}
 \end{theorem}

\section{Stability of property (R) under compact perturbations}\label{sec 4}

Let $T\in B(H)$. $T$ is said to maintain stability of property (R) under compact perturbations if $T+K\in(R)$ for any $K\in K(H)$, where $K(H)$ denotes the set of all compact operators acting on $H$.
In this section, we characterize the stability of property (R) under compact perturbations via the spectral properties of CW, and further reformulate this stability in terms of the (R)-stability of the system $\mathcal{F}_T$.

\begin{theorem} 
\label{thm4.1}
Let $T\in B(H)$. Then $T$ maintains stability of property (R) under compact perturbations if and only if $\sigma_{w}(T)=\partial\rho_{_{CW}}(T)\cup\hbox{acc}\sigma_{SF}(T)\cup\hbox{acc}\rho_{e}^{+}(T)\cup\pi_{ind}^+(T)$, where $\pi^+(T)=\{\lambda\in{\Bbb C}:n(T-\lambda )>d(T-\lambda )\}$ and $\rho_{e}^{+}(T)=\rho_e(T)\cap\pi^+(T)$.
 \end{theorem}

\begin{proof}  
 For necessity, suppose $\lambda_0\notin \partial\rho_{_{CW}}(T)\cup\hbox{acc}\sigma_{SF}(T)\cup\hbox{acc}\rho_{e}^{+}(T)\cup\pi^+(T)$. Then $n(T-\lambda_0 )\leq d(T-\lambda_0 )$.
In addition, claim $\lambda_0\in\hbox{int}\rho_{_{CW}}(T)$.
Indeed, if $\lambda_0\notin\hbox{int}\rho_{_{CW}}(T)$, then $\lambda_0\in\sigma_{CW}(T)$. Meanwhile, note that $ n(T-\lambda_0 )\leq d(T-\lambda_0 )$, thus $T-\lambda_0 \in SF_+(H)$ and $d(T-\lambda_0 )=\infty$, which means $\lambda_0\in\rho_{SF}^{-}(T)$. This contradicts ``$T$ maintains stability of property (R) under compact perturbations".
The assertion is proved.
Moreover, combining with $\lambda_0\in\hbox{int}\rho_{_{CW}}(T)$ and $\lambda_0\notin\hbox{acc}\sigma_{SF}(T)$, there exists a $\varepsilon>0$ such that $T-\lambda $ is Fredholm whenever $0<\lvert\lambda-\lambda_0\rvert<\varepsilon$. Also, note that $\lambda_0\notin\hbox{acc}\rho_{e}^{+}(T)$, thus $T-\lambda $ is Fredholm with $\hbox{ind}(T-\lambda )\leq0$ when $\lvert\lambda-\lambda_0\rvert$ is small enough. 
Therefore, $T-\lambda $ is Weyl given that $T$ maintains stability of property (R) under compact perturbations.
 It deduces $\lambda_0\in\hbox{iso}\sigma_{w}(T)\cup\rho_{w}(T)$. According to \cite[Theorem 1.1]{ba}, $\hbox{iso}\sigma_{w}(T)=\emptyset$, thus $\lambda_0\in\rho_{w}(T)$. In other words, $\lambda_0\notin\sigma_{w}(T)$. Hence 
$\sigma_{w}(T)\subseteq\partial\rho_{_{CW}}(T)\cup\hbox{acc}\sigma_{SF}(T)\cup\hbox{acc}\rho_{e}^{+}(T)\cup\pi^+(T)$. The reverse is clear.

For sufficiency, since $\hbox{iso}\sigma_{w}(T)\subseteq\hbox{int}\rho_{_{CW}}(T)$, $\rho_{SF}^{-}(T)\subseteq\hbox{int}\rho_{_{CW}}(T)\cup\sigma_{CW}(T)$, $[\hbox{iso}\sigma_{w}(T)\cup\rho_{SF}^{-}(T)]\cap[\partial\rho_{_{CW}}(T)\cup\hbox{acc}\sigma_{SF}(T)\cup\hbox{acc}\rho_{e}^{+}(T)\cup\pi^+(T)]=\emptyset$.
Thus $\hbox{iso}\sigma_{w}(T)\cup\rho_{SF}^{-}(T)=[\hbox{iso}\sigma_{w}(T)\cup\rho_{SF}^{-}(T)]\cap\sigma_{w}(T)=\emptyset$, which shows $\hbox{iso}\sigma_{w}(T)=\rho_{SF}^{-}(T)=\emptyset$. It infers from \cite[Theorem 1.1]{ba} that $T$ maintains stability of property (R) under compact perturbations.
\end{proof}

Suppose $T$ maintains stability of property (R) under compact perturbations, it follows from Theorem 4.1 that $\sigma_{w}(T)\subseteq
\overline{\sigma_{CW}(T)}\cup\hbox{acc}\sigma_{SF}(T)\cup\hbox{acc}\rho_{e}^{+}(T)$. However, the converse does not hold. For instance, let
$T\in B(\ell^2)$ be defined by $T(x_1, x_2, x_3, \cdots)=(0, x_1, 0, x_2, \cdots)$. Then $\sigma_{w}(T)=\overline{\sigma_{CW}(T)}\cup\hbox{acc}\sigma_{SF}(T)\cup\hbox{acc}\rho_{e}^{+}(T)=\Bbb D$, which shows $\hbox{iso}\sigma_{w}(T)=\emptyset$.
 Whereas, $ \rho_{SF}^{-}(T)=\{\lambda\in{\Bbb C}:0<|\lambda|<1\}\neq\emptyset$. Thus $T$ does not maintain stability of property $(R)$. To this end, let's consider the conditions such that the sufficiency holds.

\begin{corollary} 
\label{co4.3} 
Let $T\in B(H)$. Then $T$ maintains stability of property (R) under compact perturbations if and only if $\sigma_{w}(T)=\overline{\sigma_{CW}(T)}\cup\hbox{acc}\sigma_{SF}(T)\cup\hbox{acc}\rho_{e}^{+}(T)$ and $\sigma_{CW}(T)=\rho_{SF_-}(T)\cap\sigma_{SF_+}(T)$.
 \end{corollary}
\begin{proof}  
Suppose $T$ maintains stability of property (R) under compact perturbations. Based on the definition of $\sigma_{CW}(T)$, $\sigma_{CW}(T)=\rho_{SF_-}(T)\cap\sigma_{SF_+}(T)$. 
Seeing that 
$\partial\rho_{_{CW}}(T)\cup\pi_{ind}^+(T)\subseteq\overline{\sigma_{CW}(T)}\cup\hbox{acc}\rho_{e}^{+}(T)$, 
it follows from Theorem 4.1 that $\sigma_{w}(T)\subseteq\overline{\sigma_{CW}(T)}\cup\hbox{acc}\sigma_{SF}(T)\cup\hbox{acc}\rho_{e}^{+}(T)$, and the reverse is obvious.

For sufficiency, since $\hbox{iso}\sigma_{w}(T)\subseteq \hbox{int}\rho_{_{CW}}(T)$, it infers $\hbox{iso}\sigma_{w}(T)\cap[\overline{\sigma_{CW}(T)}\cup\hbox{acc}\sigma_{SF}(T)\cup\hbox{acc}\rho_{e}^{+}(T)]=\emptyset$. Thus $\hbox{iso}\sigma_{w}(T)\cap\sigma_{w}(T)=\hbox{iso}\sigma_{w}(T)=\emptyset$.
In what follows, suppose $\lambda_0\in\rho_{SF}^{-}(T)$, then $\lambda_0\in\rho_{e}(T)$. Otherwise, $T-\lambda_0 \in SF_+(H)$ and $d(T-\lambda_0 )=\infty$, which deduces $\lambda_0\in\sigma_{CW}(T)$. However, $\sigma_{CW}(T)=\rho_{SF_-}(T)\cap\sigma_{SF_+}(T)$, which is a contradiction.
Thus $\lambda_0\in\hbox{int}\rho_{_{CW}}(T)$, furthermore, $\lambda_0\notin\overline{\sigma_{CW}(T)}$. Meanwhile, $\lambda_0\notin\hbox{acc}\sigma_{SF}(T)\cup\hbox{acc}\rho_{e}^{+}(T)$. Thus $T-\lambda_0 $ is Weyl, which contradicts ``$\lambda_0\in\rho_{SF}^{-}(T)$". Hence $\rho_{SF}^{-}(T)=\emptyset$. From  \cite[Theorem 1.1]{ba}, $T$ maintains stability of property (R) under compact perturbations.
\end{proof}  

In what follows, we begin with some preparations for characterizing the stability of property (R) under compact perturbations in terms of the (R)-stability of $\mathcal{F}_T$.

\begin{lemma} 
\label{lem4.4}
Let $T\in B(H)$. If $T$ maintains stability of property (R) under compact perturbations, then $\hbox{acc}\sigma_{w}(T)=[\rho_{_{CW}}(T)\cap\pi_{\infty}(T)]\cup[\hbox{acc}\sigma(T)\cap\sigma_{d}(T)]\cup[\hbox{acc}\sigma_{a}(T)\cap\pi^+(T)]$.\ \ \ \ \ \ \ \ \ \ \ \ \ \ \ \ \ \ \ \ \ \ \ \ \ \ \ \ \ \ \ \ \ \ \ \ \ \  \ \ \ \ \ \ \ \ \ \ \ \  \ \ \ \ \ \ \ \ \ \ \ \ \ \ \ \ \ \ \ \ \ \ \ \ \ \ \ \ \ \ \ \ \ \ \ \ \ \ \ \ $(4. 1)$
 \end{lemma}

\begin{proof}  
Since $T\in(R)$, it follows from Proposition 2.7 that $\sigma_{w}(T)=[\rho_{_{CW}}(T)\cap\pi_{\infty}(T)]\cup\pi_{0}(T)\cup[\hbox{acc}\sigma(T)\cap\sigma_{d}(T)]\cup[\hbox{acc}\sigma_{a}(T)\cap\pi^{\pm}(T)]$.
Given that $T$ maintains stability of property (R) under compact perturbations, it deduces $\hbox{iso}\sigma(T)=\sigma_0(T)$, $\sigma_{w}(T)=\hbox{acc}\sigma_{w}(T)$ and $\sigma_{a}(T)=\sigma(T)$. 
  Thus $\pi_{0}(T)\subseteq\hbox{acc}\sigma(T)\cap\sigma_{d}(T)$.
 In addition, note that $\hbox{acc}\sigma_{a}(T)\cap\pi^{\pm}(T)=[\hbox{acc}\sigma_{a}(T)\cap\pi^+(T)]\cup[\hbox{acc}\sigma_{a}(T)\cap\pi^-(T)]\subseteq[\hbox{acc}\sigma_{a}(T)\cap\pi^+(T)]\cup[\hbox{acc}\sigma(T)\cap\sigma_{d}(T)]$, where $\pi^-(T)=\{\lambda\in{\Bbb C}:n(T-\lambda )< d(T-\lambda )\}$. So, $\hbox{acc}\sigma_{w}(T)\subseteq[\rho_{_{CW}}(T)\cap\pi_{\infty}(T)]\cup[\hbox{acc}\sigma(T)\cap\sigma_{d}(T)]\cup[\hbox{acc}\sigma_{a}(T)\cap\pi^+(T)]$. The anti-inclusion is clear. 
\end{proof}

 \begin{remark}
\begin{enumerate}
\item If equation $(4. 1)$ holds, then $\rho_{SF}^{-}(T)=\emptyset$. However, $\hbox{iso}\sigma_{w}(T)$ is not necessarily empty. For instance, let $T\in B(\ell^2)$ be defined by $T(x_1, x_2, x_3, \cdots)=(0, 0, \frac{x_2}{2},\frac{x_3}{3}, \cdots)$.
Then $(4. 1)$ holds, whereas $\hbox{iso}\sigma_{w}(T)=\{0\}$. Thus	$T$ does not maintain stability of property (R) under compact perturbations.						

\item If $T$ maintains stability of property (R) under compact perturbations, then
$$\begin{array}{rcl}
\hbox{acc}\sigma_{w}(T)&=&[\rho_{_{CW}}(T)\cap\pi_{\infty}(T)]\cup[\hbox{acc}\sigma(T)\cap\sigma_{d}(T)]\cup\pi_0^a(T)\\&\cup&[\hbox{acc}\sigma_{a}(T)\cap\pi^+(T)],\ \ \ \ \ \ \ \ \ \ \ \ \ \ \ \ \ \ \  \ \ \ \ \ \ \ \ \ \ \ \ \ \ \ \ \ \ \ \ \ \ \ \ \ \ (4. 2) 
\end{array}$$
where $\pi_0^a(T)=\{\lambda\in\sigma_{a}(T):n(T-\lambda )=0\}$.
However, the reverse doesn't hold. Indeed,
assume that equation $(4. 2)$ holds. Then $\rho_{SF}^{-}(T)=\emptyset$. However, $\hbox{iso}\sigma_{w}(T)$ may be non-empty, as shown in the example mentioned in Remark 4.4 (1). Thus $T$ does not maintain stability of property (R) under compact perturbations. 
\end{enumerate}
\end{remark}

 To characterize the stability of property (R) under compact perturbations in terms of the (R)-stability of $\mathcal{F}_T$, we clarify the relationship between these two notions of stability.

Suppose $T$ maintains stability of property (R) under compact perturbations, then $\mathcal{F}_T$ is (R)-stable.
 Indeed, it follows from \cite[Theorem 1.1]{ba} that $\sigma_{ab}(T)=\sigma_{b}(T)$ and $\hbox{iso}\sigma(T)=\sigma_0(T)$. In addition, due to Lemma 3.3, the system $\mathcal{F}_T$ is (R)-stable.
Conversely, suppose that $\mathcal{F}_T$ is (R)-stable, the following examples show that $T$ may not maintain stability of property (R) under compact perturbations. 
\begin{example}
\label{ex4.6}
Let $A, B\in B(\ell^2)$ be defined by $A(x_1, x_2, x_3, \cdots)=(0, x_1, x_2, \cdots)$, $B(x_1, x_2, x_3, \cdots)=(0, x_1, \frac{x_2}{2},\frac{x_3}{3}, \cdots)$. Let $T_1=A\oplus(B+3)$.
It can be checked that $\sigma_0(T_1)=\emptyset$, $\sigma_{a}(T_1)=\sigma_{ab}(T_1)=\Bbb T_1\cup\{3\}$, $\pi_{00}(T_1)=\emptyset$, which means $T_1\in(R)$. From Lemma 3.3, $\mathcal{F}_{T_1}$ is (R)-stable.
Whereas, 
$\hbox{iso}\sigma_{w}(T_1)=\{3\}\neq\emptyset$, $\rho_{SF}^{-}(T_1)=\{\lambda\in{\Bbb C}:|\lambda|<1\}\neq\emptyset$. Thus $T_1$ does not maintain stability of property (R) under compact perturbations on account of \cite[Theorem 1.1]{ba}. Besides that, $\hbox{acc}\sigma_{w}(T_1)=\Bbb D$ and $[\rho_{_{CW}}(T_1)\cap\pi_{\infty}(T_1)\}]\cup[\hbox{acc}\sigma(T_1)\cap\sigma_{d}(T_1)]\cup[\hbox{acc}\sigma_{a}(T_1)\cap\pi^+(T_1)]\cup\pi_0^a(T_1)=\Bbb T\cup\{3\}$.
\end{example}

\begin{example}
\label{ex4.7}
Let $A, B, C, D\in B(\ell^2)$ be defined by $A(x_1, x_2, x_3, \cdots)=(0, 0, x_1, x_2, \cdots)$, $B(x_1, x_2, x_3, \cdots)=(0, x_1, 0, \frac{x_3}{3}, 0, \frac{x_5}{5}, \cdots)$, $C(x_1, x_2, x_3, \cdots)=(x_2, x_3, x_4, \cdots)$, $D(x_1, x_2, x_3, \cdots)=(x_1, 0, 0, \cdots)$. Let $T_2=A\oplus (B+3)\oplus C\oplus (D+1)$.
Now, $\sigma_0(T_2)=\{2\}\neq\emptyset$. Furthermore, $\sigma_{ab}(T_2)=\sigma_{b}(T_2)=\Bbb D\cup\{3\}$, $\hbox{iso}\sigma(T_2)\backslash \pi_{\infty}(T_2)=\sigma_0(T_2)=\{2\}$. Combining with Lemma 3.3, $\mathcal{F}_{T_2}$ is (R)-stable.
However, note that $\hbox{iso}\sigma_{w}(T_2)=\{3\}\neq\emptyset$, $\rho_{SF}^{-}(T_2)=\Bbb D\backslash\Bbb T\neq\emptyset$. It infers from \cite[Theorem 1.1]{ba} that $T_2$ does not maintain stability of property (R) under compact perturbations. Besides that,  $\hbox{acc}\sigma_{w}(T_2)=\Bbb D$ and $[\rho_{_{CW}}(T_2)\cap\pi_{\infty}(T_2)\}]\cup[\hbox{acc}\sigma(T_2)\cap\sigma_{d}(T_2)]\cup[\hbox{acc}\sigma_{a}(T_2)\cap\pi^+(T_2)]\cup\pi_0^a(T_2)=\Bbb T\cup\{3\}$.
\end{example}

One can see that neither the $T_1$  in Example 4.5 nor the $T_2$ in Example 4.6 satisfies the equality $(4. 2)$. However, Remark 4.4 (2) shows that the $(4. 2)$ is closely related to the stability of property (R) under compact perturbations.
Motivated by the preceding analysis, 
the stability of property (R) under compact perturbations is established as detailed below.

\begin{theorem} 
\label{thm4.8}
Let $T\in B(H)$. Then $T$ maintains stability of property (R) under compact perturbations if and only if the following conditions hold:
\begin{enumerate} 
\item The system $\mathcal{F}_T$ is (R)-stable;
\item $\hbox{acc}\sigma_{w}(T)=[\rho_{_{CW}}(T)\cap\pi_{\infty}(T)\}]\cup[\hbox{acc}\sigma(T)\cap\sigma_{d}(T)]\cup[\hbox{acc}\sigma_{a}(T)\cap\pi^+(T)]\cup\pi_0^a(T).\quad \quad \quad \quad \quad \quad \quad \quad \quad \quad \quad \quad \quad \quad \quad \quad \quad \quad \quad \quad \quad \quad \quad \quad \quad \quad \quad \quad\quad   (4. 2)$
\end{enumerate}
 \end{theorem}

\begin{proof}
By Remark 4.4 (2), it suffices to establish the sufficiency by showing that $\hbox{iso}\sigma_{w}(T)=\emptyset$. Indeed, assume on the contrary that 
$\lambda_0\in\hbox{iso}\sigma_{w}(T)$,  then $\lambda_0\notin\hbox{acc}\sigma_w(T)$, and thus $\lambda_0\notin[\rho_{_{CW}}(T)\cap\pi_{\infty}(T)]\cup[\hbox{acc}\sigma(T)\cap\sigma_{d}(T)]\cup[\hbox{acc}\sigma_{a}(T)\cap\pi^+(T)]\cup\pi_{0}^a(T)$. 
Note that $\hbox{iso}\sigma_{w}(T)\subseteq\rho_{_{CW}}(T)$, thus $\lambda_0\notin\pi_{\infty}(T)$, i.e., $n(T-\lambda_0 )<\infty$. 
Based on the fact that $\hbox{iso}\sigma_{w}(T)\subseteq\sigma_{a}(T)$, one can see $\lambda_0\in\sigma_{a}(T)$. Considering that $\lambda_0\notin\pi_{0}^a(T)$,  $n(T-\lambda_0 )>0$. So, $0<n(T-\lambda_0 )<\infty$. 
Furthermore, it infers from $\lambda_0\notin\hbox{acc}\sigma(T)\cap\sigma_{d}(T)$ that either $\lambda_0\in\pi_{00}(T)$ or $ T-\lambda_0 \in SF_+(H)$ with $ n(T-\lambda_0 )>0$. 
If $\lambda_0\in\pi_{00}(T)$, it infers from $T\in(R)$ that $ T-\lambda_0 $ is Browder, which contradicts ``$\lambda_0\in\hbox{iso}\sigma_{w}(T)$". 
Thus, $T-\lambda_0 \in SF_+(H)$ with $ n(T-\lambda_0 )>0$. Besides that, note that $\lambda_0\notin\hbox{acc}\sigma_{a}(T)\cap\pi^+(T)$, then $\lambda_0\in\sigma_{a}(T)\backslash\sigma_{ab}(T)$ or $\hbox{ind}(T-\lambda_0)\leq0$.
 If $\lambda_0\in\sigma_{a}(T)\backslash\sigma_{ab}(T)$, it deduces from $T\in(R)$ that $T-\lambda_0 $ is Browder, which contradicts ``$\lambda_0\in\hbox{iso}\sigma_{w}(T)$". Thus $ T-\lambda_0 \in SF_+(H)$ with $ n(T-\lambda_0 )>0$, and $\hbox{ind}(T-\lambda_0)\leq0$. From Remark 4.4 (2), $\rho_{SF}^{-}(T)=\emptyset$, which means $T-\lambda_0$ is Weyl. This again contradicts ``$\lambda_0\in\hbox{iso}\sigma_{w}(T)$".
To sum up, $\hbox{iso}\sigma_{w}(T)=\emptyset$, hence $T$ maintains stability of property (R) under compact perturbations.
\end{proof}

Combining with Remark 4.4 (2), the part ``$\pi_{0}^a(T)$" in the decomposition of $\hbox{acc}\sigma_{w}(T)$ is essential in Theorem 4.7.
Indeed, let $T\in B(\ell^2)$ be defined by $T(x_1, x_2, x_3, \cdots)=(0, x_1, \frac{x_2}{2},\frac{x_3}{3}, \cdots)$. Then $\sigma_{0}(T)=\emptyset$, $\sigma_{a}(T)=\sigma_{ab}(T)=\{0\}$ and $\pi_{00}(T)=\emptyset$. Thus $T\in(R)$, which shows $\mathcal{F}_T$ is (R)-stable.
Besides that, $\hbox{acc}\sigma_{w}(T)=\emptyset=[\rho_{_{CW}}(T)\cap\pi_{\infty}(T)]\cup[\hbox{acc}\sigma(T)\cap\sigma_{d}(T)]\cup[\hbox{acc}\sigma_{a}(T)\cap\pi^+(T)]$. However, $\hbox{iso}\sigma_{w}(T)=\{0\}\neq\emptyset$, hence $T$ does not maintain stability of property (R) under compact perturbations. 
To omit the part ``$\pi_{0}^a(T)$", we now focus on the stability of property (R) under compact perturbations under the condition that $\sigma_0(T)\neq\emptyset$.

\begin{corollary} 
\label{co4.9}
Let $T\in B(H)$. If $\sigma_0(T)\neq\emptyset$, then $T$ maintains stability of property (R) under compact perturbations if and only if the following conditions hold:
\begin{enumerate} 
\item The system $\mathcal{F}_T$ is (R)-stable;
\item $\hbox{acc}\sigma_{w}(T)=[\rho_{_{CW}}(T)\cap\pi_{\infty}(T)]\cup[\hbox{acc}\sigma(T)\cap\sigma_{d}(T)]\cup[\hbox{acc}\sigma_{a}(T)\cap\pi^+(T)]$. \ \ \ $(4. 1)$
\end{enumerate}
 \end{corollary}

\begin{proof}
It suffices to show the sufficiency. Considering that $\sigma_{0}(T)\neq\emptyset$ and $\mathcal{F}_T$ is (R)-stable, one can obtain that $\sigma_{ab}(T)=\sigma_{b}(T)$ and $\hbox{iso}\sigma(T)\backslash \pi_{\infty}(T)=\sigma_{0}(T)$, thus $\sigma_{a}(T)=\sigma(T)$. Assume on the contrary that $\lambda_0\in\hbox{iso}\sigma_{w}(T)$, it deduces from the proof of Theorem 4.7 that either $\lambda_0\in\hbox{iso}\sigma(T)\backslash \pi_{\infty}(T)$ or $ T-\lambda_0\in SF_+(H)$ with $n(T-\lambda_0)>0$. Combining with $\hbox{iso}\sigma(T)\backslash \pi_{\infty}(T)=\sigma_{0}(T)$ and $\lambda_0\in\hbox{iso}\sigma_{w}(T)$, $ T-\lambda_0\in SF_+(H)$ with $n(T-\lambda_0)>0$.  In what follows, analogous to the proof of Theorem 4.7, one can get that $\hbox{iso}\sigma_{w}(T)=\emptyset$. In addition, it infers from the equation $(4. 1)$ that $\rho_{SF}^{-}(T)=\emptyset$. Therefore, $T$ maintains stability of property (R) under compact perturbations.
\end{proof}

\section*{Acknowledgements}
The authors would like to express their sincere gratitude to Professor Xu Houbao for providing the research platform and for his support.
\section*{Data availability}

All data included in this study are available upon request by contacting with the corresponding author.

\section*{Declarations}

There is no potential conflict of interest between the authors of this paper.

\end{document}